\documentclass[12pt]{article}

\usepackage{latexsym}
\usepackage{graphicx}
\usepackage{color}

\newcommand{\bb}{\mbox{\boldmath$b$}}
\newcommand{\be}{\mbox{\boldmath$e$}}

\newcommand{\bp}{\mbox{\boldmath$p$}}
\newcommand{\bq}{\mbox{\boldmath$q$}}
\newcommand{\br}{\mbox{\boldmath$r$}}
\newcommand{\bR}{{\bf R}}
\newcommand{\bs}{\mbox{\boldmath$s$}}

\newcommand{\bv}{\mbox{\boldmath$v$}}
\newcommand{\bx}{\mbox{\boldmath$x$}}
\newcommand{\by}{\mbox{\boldmath$y$}}
\newcommand{\bze}{{\bf 0}}

\newcommand{\nul}{{\cal N}}

\newcommand{\ran}{{\cal R}}
\newcommand{\rank}{\mbox{rank}}
\newcommand{\rmn}{\bR^{m \times n}}
\newcommand{\rn}{{\bf R}^n}

\newcommand{\rnn}{{\bf R}^{n \times n}}

\newcommand{\trans}{{\mbox{\scriptsize T}}}

\makeatletter
 
 \@addtoreset{equation}{section}
\makeatother
\makeatletter
 
 \@addtoreset{figure}{section}
\makeatother
\newfont{\bg}{cmr10 scaled\magstep4}

\newcommand{\bigzerou}{\smash{\lower1.7ex\hbox{\bg 0}}}
\begin{document}
\title{Convergence of the Conjugate Gradient Method \\ on Singular Systems\footnote{The contents of this paper was presented as a part of \cite{HY}}}
\author{Ken Hayami\footnote{National Institute of Informatics, 2-1-2, 
Hitotsubashi, Chiyoda-ku, Tokyo 101-8430, Japan, and 
SOKENDAI (The Graduate University for Advanced Studies), 
e-mail: hayami@nii.ac.jp} }
\date{}
\maketitle
\begin{abstract}
We analyze the convergence of the Conjugate Gradient (CG) method in exact arithmetic, when the coefficient matrix $A$ is symmetric positive semidefinite and the system is consistent. To do so, we diagonalize $A$ and decompose the algorithm into the range and the null space components of $A$.
Further, we apply the analysis to the CGLS and CGNE (CG Normal Error) methods for rank-deficient least squares problems.
\end{abstract}

\noindent
{\bf Key words:}\hspace{2mm}conjugate gradient method, singular systems, least squares problems\\
\noindent
{\bf Mathematical Sub Classification Number:}\hspace{2mm}65F10

\section{Introduction}
\label{intro}

In exact arithmetic, the Conjugate Gradient (CG) method \cite{HeS} converges 
to the true solution within $n$ iterations for a system of linear equations 
\begin{equation}
   A \bx = \bb, \hspace{5mm}A \in \rnn ,
 \label{CG}
\end{equation}
if $A$ is symmetric positive definite.

When $A$ is symmetric positive semidefinite, the CG method still converges to a solution, 
provided $ \bb \in \ran (A) $. (This fact is stated, for instance, in Kaasschieter \cite{K}, p.266.) In this short paper, we will show this explicitly by diagonalizing $A$ and decomposing the CG method into the 
$\ran (A)$ component and the $\nul (A)=\ran (A)^\perp$ component, where 
$\ran (A)$ is the range space of $A$, and $\nul (A)$ is the null space of $A$, 
and $S^\perp$ is the orthogonal complement of a subspace $S$. This follows the 
analysis of the Conjugate Residual method  \cite{A99, H00, H01, H03} and GMRES (Generalized Minimal Residual) and GCR (Generalized Conjugate Residual) methods in \cite{HaS}, 
but the analysis is much simpler in the symmetric case.

Using this fact, we will also show that the CGLS and CGNE (CG Normal Error) methods converge 
for least squares problems 
\[ \min_{\bx \in \rn} || \bb - A \bx ||_2, \hspace{5mm}A \in \rmn , \]
even when $A$ is rank-deficient.
\section{The CG algorithm}

The CG algorithm \cite{HeS} for solving (\ref{CG}) is given as follows, where $\bx_0$ is the initial approximate solution.

\noindent
\hspace{1.5cm}{\bf The CG algorithm}\\ \\
$\br_0=\bb-A\bx_0\\
\bp_0=\br_0$\\
For $i=0,1,\ldots$, until converegence, Do\\
\hspace{1cm}$\alpha_i = \frac{(\br_i,\br_i)}{ {\displaystyle (A \bp_i, \bp_i)} }\\
\hspace{1cm}\bx_{i+1}=\bx_i + \alpha_i \bp_i\\
\hspace{1cm}\br_{i+1}=\br_i-\alpha_i A\bp_i\\
\hspace{1cm}\beta_i=\frac{(\br_{i+1},\br_{i+1})}{(\br_i,\br_i)}$\\
\hspace{1cm}$\bp_{i+1}=\br_{i+1}+\beta_i\bp_i$\\ 
End Do \\

\section{Decomposition of the CG algorithm into the range and null space components}
In the following, we assume exact arithmetic. 
When $A$ is symmetric positive definite, the CG method is known to converge to 
the true solution within $n$ iterations \cite{HeS}. When $A$ is symmetric 
positive semidefinite, we will analyze the method as follows.

Let 
\[ Q^\trans A Q= \Lambda=
   \left[ \begin{array}{cccc}
                         \Lambda_r &   &        &   \\
                                   & 0 &        &   \\  
                                   &   & \ddots &   \\
                                   &   &        & 0 
          \end{array}
    \right],
\]
where
\[ \Lambda_r = \left[ \begin{array}{ccc}
                      \lambda_1 &        &            \\
                                & \ddots &            \\
                                &        & \lambda_r 
                      \end{array}
               \right],\hspace{1cm}
\lambda_1 \ge \lambda_2 \cdots \ge \lambda_r > 0, 
\hspace{1cm} r = \rank A = \dim \ran (A), \]
and
\[ Q^\trans Q = I_n, \hspace{1cm} Q=\left[ Q_1, Q_2 \right], \]
\[ Q_1=\left[ \bq_1, \cdots, \bq_r \right], \hspace{1cm}
   Q_2=\left[ \bq_{r+1}, \cdots, \bq_n \right]. \]

Here, $\bq_1, \cdots, \bq_r$ and $\bq_{r+1}, \cdots, \bq_n
$ are orthonormal basis vectors of $\ran(A)$ and $\ran(A)^\perp = \nul(A)$, respectively.
   
Following \cite{HaS}, we will decompose the CG algorithm into the $\ran(A)$ component and the $\ran(A)^\perp = \nul(A)$ component using the following transformation:
\[ \tilde{\bv} = Q^\trans \bv = [ Q_1 , Q_2 ]^\trans \bv
               = \left[ \begin{array}{c}
                          {Q_1}^\trans \bv \\
                          {Q_2}^\trans \bv    
                        \end{array}
                 \right]
               = \left[ \begin{array}{c}
                           \bv^1 \\
                           \bv^2    
                        \end{array}
                  \right]
\]
or
\[ \bv= Q \tilde{\bv}= [ Q_1 , Q_2 ] \left[ \begin{array}{c}
                                              \bv^1 \\
                                              \bv^2    
                                            \end{array}
                                     \right]
       = Q_1 \bv^1 + Q_2 \bv^2 .
\]
Here, $Q_1 \bv^1$ and $Q_2 \bv^2 $ correspond to the $\ran(A)$ and $\nul(A)$ component of $\bv$, respectively.

Noting that
\[ 
\begin{array}{lll}
     \br & = & \bb - A \bx, \\ \\
     Q^\trans \br & = & Q^\trans \bb - Q^\trans A Q ( Q^\trans \bx ), \\ \\   
     \left[ \begin{array}{c}
                           \br^1 \\
                           \br^2    
                        \end{array}
                  \right] &
     =  & \left[ \begin{array}{c}
                           \bb^1 \\
                           \bb^2    
                        \end{array}
                  \right]
        - \left[
              \begin{array}{cc}
               \Lambda_r &  0\\
                0        &  0   
              \end{array}
             \right]
          \left[   
            \begin{array}{c}
                           \bx^1 \\
                           \bx^2    
                        \end{array}
          \right], \\ \\
 \br^1 & = & \bb^1 - \Lambda_r \bx^1, \\ 
 \br^2 & = & \bb^2, \\ \\
(\br,\br) &=&(\br^1,\br^1)+(\bb^2,\bb^2),\\ \\ 
(A\bp,\bp)&=&(\bp^1,\Lambda_r \bp^1),
\end{array}
\]
we obtain the following.\\

\hspace{3cm}{\bf Decomposed CG algorithm\hspace*{5mm}(General case)}
\[
\begin{array}{lll}
\underline{ {\cal R}(A) \, \mbox{component} } & \hspace{5mm} & 
\underline{ \nul(A) \, \mbox{component} } \\ \\
\br_0^1 = \bb^1 - \Lambda_r \bx_0^1 
& \hspace{0cm} & \br_0^2 = \bb^2 \\ 
\bp_0^1 = \br_0^1 & \hspace{0cm} & 
\bp_0^2 = \bb^2 \\ \\
\multicolumn{3}{l}{\mbox{For }i = 0,1,\ldots, \mbox{until convergence
Do} } \\ \\
\hspace{1cm}\alpha_i = \frac{ ( \br_i^1, \br_i^1) + (\bb^2,\bb^2) }
                {\displaystyle ( \bp_i^1, \Lambda_r \bp_i^1 ) } & & \\  
\hspace{1cm}\bx_{i+1}^1 = \bx_i^1 + \alpha_i \bp_i^1 & \hspace{0cm} &
\bx_{i+1}^2 = \bx_i^2 + \alpha_i \bp_i^2 \\ 
\hspace{1cm}\br_{i+1}^1 = \br_i^1 - \alpha_i {\displaystyle \Lambda_r \bp_i^1}
 & \hspace{0cm} & \br_{i+1}^2 = \bb^2 \\ 
\multicolumn{3}{l}
{\hspace{1cm}\beta_i = \frac{ ( \br_{i+1}^1, \br_{i+1}^1) +  ( \bb^2, \bb^2 ) }
                 { ( \br_i^1, \br_i^1) +  ( \bb^2, \bb^2 ) } }\\ 
\hspace{1cm}\bp_{i+1}^1 = \br_{i+1}^1 + \beta_i \bp_i^1 & \hspace{0cm} &
\bp_{i+1}^2 = \bb^2 + \beta_i \bp_i^2 \\ 
\multicolumn{3}{l}{\mbox{End Do}}
\end{array}
\]

Here, it is interesting to note that $\br_i^2=\bb^2$ for $i=1,2,\ldots$, and $\bp_i^2, \: i=1,2,\ldots$ is confined in a one-dimensional subspace parallel to $\bb^2$ within the $(n-r)$-dimensional subspace corresponding to ${\cal N}(A)$.

Further, if we assume that the system is consistent, i.e., 
$ \bb \in \ran (A)$ or $\bb^2 = \bze$, we have the following.

\hspace{2cm}
{\bf Decomposed CG algorithm\hspace*{5mm}($\bb \in \ran(A)$: $\bb^2=\bze$)}
\[
\begin{array}{lll}
\underline{ {\cal R}(A) \, \mbox{component} } & \hspace{5mm} & 
\underline{ \nul(A) \, \mbox{component} } \\ \\
\br_0^1 = \bb^1 - \Lambda_r \bx_0^1 
& \hspace{0cm} & \br_0^2 = \bze\\ 
\bp_0^1 = \br_0^1 & \hspace{0cm} & 
\bp_0^2 = \bze\\ \\
\multicolumn{3}{l}{\mbox{For } \, i = 0,1, \, \ldots , \mbox{until convergence
Do} } \\ \\
\hspace{1cm}\alpha_i = \frac{ ( \br_i^1, \br_i^1) }
                {\displaystyle ( \bp_i^1, \Lambda_r \bp_i^1 ) } & & \\  
\hspace{1cm}\bx_{i+1}^1 = \bx_i^1 + \alpha_i \bp_i^1 & \hspace{0cm} &
\bx_{i+1}^2 = \bx_0^2\\ 
\hspace{1cm}\br_{i+1}^1 = \br_i^1 - \alpha_i {\displaystyle \Lambda_r \bp_i^1}
 & \hspace{0cm} & \br_{i+1}^2 = \bze \\ 
\multicolumn{3}{l}
{\hspace{1cm}\beta_i = \frac{ ( \br_{i+1}^1, \br_{i+1}^1) }
                 { ( \br_i^1, \br_i^1) } }\\ 
\hspace{1cm}\bp_{i+1}^1 = \br_{i+1}^1 + \beta_i \bp_i^1 & \hspace{0cm} &
\bp_{i+1}^2 = \bze\\ 
\multicolumn{3}{l}{\mbox{End Do}}
\end{array}
\]

Note that in the above, the $\ran(A)$ component is equivalent to CG applied to $\Lambda_r \bx^1 = \bb^1$, where $\Lambda_r$ is symmetric positive definite. 
Hence, $ \bx_i^1 $ will converge to $ \bx_\ast^1={\Lambda_r}^{-1}\bb^1 $ within $r=\rank A$ iterations. Thus, $ \br_i^1 $ will converge to $ \bze$.
Therefore, the $\ran(A)$ component $Q_1 \bx_i^1$ of $\bx_i$ will converge to $A^\dagger \bb$, and the $\ran(A)$ component of $\br_i$ will converge to $\bze$. 

As for the $\nul(A)$ component, since $\bx_i^2=\bx_0^2$ and $\br_i^2=\bze$ for $i=0,1, \ldots$ , $Q_2 \bx_i^2 =Q_2 \bx_0^2$ and $Q_2 \br_i^2 = \bze$ for $i=0,1,\ldots$ . Thus, $\bx_i$ will converge to a solution of (\ref{CG}). 

If further, $\bx_0 \in \ran(A)$ (e.g. $\bx_0=\bze)$, then, $\bx_0^2=\bze$. 
Thus, $\bx_i^2=\bze$, \hspace{2mm}$Q_2 \bx_i=\bze$ for $i=0,1,\ldots$ .
Therefore, $\bx_i$ will converge to the minimum norm solution (pseudo inverse solution) $A^\dagger \bb$ of (\ref{CG}).

Analysis similar to the above was done in Washizawa et al.(\cite{WY}, Subsection 4.2). Here, we have further simplified the analysis by diagonalizing $A$ 

We also obtain the following error bound:
\[ ||\be_k^1||_{\Lambda_r} = ||\br_k^1||_{ {\Lambda_r}^{-1} }
   \le 2 \left\{ \frac{\sqrt{\kappa(\Lambda_r)}-1}{\sqrt{\kappa(\Lambda_r)}+1} \right\}^k
          ||\br_0^1||_{ {\Lambda_r}^{-1} } \, ,
\] 
where 
\[ \be_k^1=\bx_k^1-\bx_\ast^1 \, ,  \hspace{1cm} \kappa(\Lambda_r)=\frac{\lambda_1}{\lambda_r}.\]

Note here that
\[ ||\be^1||_{\Lambda_r}^2 = {\be^1}^\trans \Lambda_r \be^1 
    = \be^\trans A \be = ||\be||_A\]
and
\[ ||\br^1||_{ {\Lambda_r}^{-1} }^2 = {\br^1}^\trans {\Lambda_r}^{-1} \br^1 
 = \br^\trans A^\dagger \br.\]

Note also that
\[ A=Q \Lambda Q^\trans, \]
\[ A^\dagger=Q \Lambda^\dagger Q^\trans, \]
where\\
\[ \Lambda^\dagger
     =\left[ \begin{array}{cccc}
                  {\Lambda_r}^{-1} &   &        &   \\
                                   & 0 &        &   \\  
                                   &   & \ddots &   \\
                                   &   &        & 0 
             \end{array}
       \right].
\]
Here
\[ \bx_\ast=A^\dagger \bb = [ Q_1, Q_2 ] 
       \left[ \begin{array}{cccc}
                  {\Lambda_r}^{-1} &   &        &   \\
                                   & 0 &        &   \\  
                                   &   & \ddots &   \\
                                   &   &        & 0 
             \end{array}
       \right] 
       \left[ \begin{array}{c}
                 {Q_1}^\trans \\
                 {Q_2}^\trans
              \end{array}
       \right]
       \bb
      = Q_1 {\Lambda_r}^{-1} \bb^1 = Q_1 \bx_\ast^1 \, .
\]

\section{Application to rank-deficient least squares problems}
Next, we apply the above analysis to the least squares problem
\begin{equation}
 \min_{\bx \in \rn} \| \bb - A \bx \|_2 , 
 \label{lsq}
\end{equation}
where now $ A \in \bR^{m\times n}$ and  $m \ge n, $ and $A$ may not necessarily be of full rank.

Equation (\ref{lsq}) is equivalent to the normal equation (of the first kind):
\begin{equation}
A^\trans A \bx = A^\trans \bb .
\label{NR}
\end{equation}
Here,  $A^\trans A$ is square, symmetric, and positive semidefinite. Moreover, (\ref{NR}) is consistent, 
since $A^\trans \bb = A^\trans \bb \in \ran(A^\trans A)=\ran(A^\trans)$.
Thus, we may apply the Conjugate Gradient method to (\ref{NR}) with $A:=A^\trans A$ and $\bb:= A^\trans \bb$
 to obtain the following CGLS(CGNR(CG Normal Residual)) method\cite{B,S}.\\ \\
\hspace{3cm}{\bf CGLS method}\\ \\
\hspace{6mm} Choose $ \bx_0 $. \\
\hspace{6mm} $ \br_0 = \bb-A\bx_0,\hspace{3mm}\bp_0 = \bs_0 = A^\trans \br_0,
\hspace{3mm}\gamma_0 = {\| \bs_0 \|_2}^2 $ \\
\hspace{6mm} For $i=0,1,2,\ldots,$ until $\gamma_i < \epsilon$ \vspace{1mm} 
Do\\
\hspace{15mm}$ \bq_i = A \bp_i $ \\
\hspace{15mm}$ \alpha_i = {\gamma_i / \| \bq_i \|_2}^2 $ \\
\hspace{15mm}$ \bx_{i+1} =\bx_i + \alpha_i \bp_i $ \\
\hspace{15mm}$ \br_{i+1} = \br_i - \alpha_i \bq_i $ \\
\hspace{15mm}$ \bs_{i+1} = A^\trans \br_{i+1} $ \\
\hspace{15mm}$ \gamma_{i+1} = { \| \bs_{i+1} \|_2}^2 $ \\
\hspace{15mm}$ \beta_i = \gamma_{i+1} / \gamma_i $ \\
\hspace{15mm}$ \bp_{i+1} = \bs_{i+1} + \beta_i \bp_i $
\vspace{1mm}\\
\hspace{6mm} End Do\\

In the above algorithm, the $\ran(A^\trans)$ component of the solution $Q_1 \bx_k^1 $ converges to 
$ (A^\trans A)^\dagger A^\trans \bb = A^\dagger \bb$, and $A^\trans \br_k$ converges to $\bze$.

The $\nul(A^\trans)$ component of the solution remains $\bx_k^2 = \bx_0^2$.

Thus, if  $\bx_0 \in \ran(A^\trans)$ (e.g. $\bx_0=\bze)$, $\bx_k^2=\bze$, so that the solution $\bx_k$ converges to the minimum norm solution: $ \bx_\ast=A^\dagger \bb $.

This was shown by Kammerer and Nashed \cite{KN}(Theorem 5.1) in a Hilbert space setting. Here we gave a simpler derivation using matrix analysis.

We also obtain the following error bound:

\[ \| \br_k |_{\ran(A)} \|_2 = \be_k^\trans A^\trans A \be_k
      \le 2 \left( \frac{\sigma_1 -\sigma_r}{\sigma_1 +\sigma_r} \right)^k
       \| \br_0 |_{\ran(A)} \|_2, \]
where $\sigma_1 \ge \cdots \ge \sigma_r > 0 $ are the  nonzero singular values of $A$. \\

Similar analysis can be performed for the under-determined least squares problem:
\[ \min_{\{\bx | \displaystyle {A \bx=\bb } \}} || \bx ||_2 \]
with $m \le n$, which is equivalent to the normal equation of the second kind: 
\[ A A^\trans \by = \bb,\hspace{5mm}\bx=A^\trans \by \,  . \]
 Applying  CG to this equation with $A: AA^\trans$ gives the CGNE method \cite{S}.
 Assume $ \bb \in \ran (AA^\trans)=\ran (A)$. Further assume that the initial approximate solution satisfies 
$\by_0 \in \ran (A)$, which is equivalent to $\bx_0 \in \ran(A^\trans A)=\ran(A^\trans)$ (e.g. $\by_0=\bze$, $\bx_0=\bze$).
Then, $\by_k$ will converge to $(AA^\trans)^\dagger \bb$, and the solution $\bx_k$ will converge to $A^\trans  (AA^\trans)^\dagger \bb = A^\dagger \bb$, the min-norm solution.

 We obtain the following error bound:
\[ \br_k^\trans (AA^\trans)^\dagger \br_k \le 
     2 \left( \frac{\sigma_1 -\sigma_r}{\sigma_1 +\sigma_r} \right)^k 
      \br_0^\trans (AA^\trans)^\dagger \br_0 . \]

\section{Concluding remark}
In this paper, we analyzed the convergence of the CG method in exact arithmetic, when the coefficient matrix $A$ is symmetric positive semidefinite and the system is consistent. To do so, we diagonalized $A$ and decomposed the algorithm into the $\ran(A)$ component and its orthogonal complement space $\nul(A)=\ran(A)^\perp$ component. Further, we applied the analysis to the CGLS and CGNE methods for rank-deficient least squares problems.

\section{Acknowledgenets}
This work was supported in part by JSPS KAKENHI  Grant Number 15K04768.

\end{document}